\documentclass{amsart}

\DeclareMathOperator{\Conv}{Conv}
\DeclareMathOperator{\Lie}{Lie}
\DeclareMathOperator{\val}{val}
\DeclareMathOperator{\Hom}{Hom}

\usepackage{amssymb}

\numberwithin{equation}{section}

\newtheorem{thm}{Theorem}[section]

\newtheorem{lem}[thm]{Lemma}
\newtheorem{prop}[thm]{Proposition}

\theoremstyle{definition}

\newtheorem{question}[thm]{Question}

\begin{document}
\title{On the Hodge-Newton decomposition for split groups}

\author{Robert E. Kottwitz
\address {Department of Mathematics, University of Chicago, 5734 University
Avenue, Chicago, Illinois 60637}
}

\email{kottwitz@math.uchicago.edu}
\thanks{Partially supported by NSF Grant DMS-0071971}


\maketitle

\section{Introduction} The main purpose of this paper is to prove a group-theoretic
generalization of a theorem of Katz on isocrystals. Along the way we reprove  the
group-theoretic generalization of Mazur's inequality for isocrystals 
 due to Rapoport-Richartz, and generalize from split groups to unramified groups
a result from \cite{kottwitz-rapoport02} which determines when the affine
Deligne-Lusztig subset $X^G_\mu(b)$ of $G(L)/G(\mathfrak o_L)$ is non-empty.

 Let $F$ be a finite extension of~$\mathbf Q_p$ with uniformizing element
$\varpi$. We write
$L$ for the completion of the maximal unramified extension of~$F$ in some
algebraic closure~$\bar F$ of~$F$. We write $\sigma$ for the Frobenius
automorphism of~$L$ over~$F$, and we write 
$\mathfrak o$ (respectively, $\mathfrak o_L$) for the valuation ring of~$F$
(respectively,
$L$).

Let $G$ be a split connected reductive group over $\mathfrak o$ and let $A$ be a
split maximal torus of~$G$ over~$\mathfrak o$.  Fix a Borel subgroup $B=AU$
containing $A$ with unipotent radical~$U$, as well as a parabolic subgroup~$P$
of~$G$ containing~$B$. Write
$P=MN$, where $M$ is the unique Levi subgroup of~$P$ containing $A$ and $N$ is the
unipotent radical of~$P$.

We write $X_G$ for the 	quotient of $X_*(A)$ by the coroot lattice for $G$, and we
define  a homomorphism $w_G:G(L) \to X_G$ as follows. For $g\in G(L)$  we define
$r_B(g)
\in X_*(A)$ to be the unique element
$\mu \in X_*(A)$ such that $g \in G(\mathfrak o_L) \cdot \mu(\varpi) \cdot U(L)$,
and we define
$w_G(g)$ to be the image of $r_B(g)$ under the canonical surjection from $X_*(A)$
to $X_G$. (This definition of $w_G$ suffices for the purposes of this paper; see
\S7 of
\cite{kottwitz97} for a definition that applies to groups $G$ that are not split
over $L$.)

Applying the construction above to $M$ rather than $G$, we obtain $X_M$, the
quotient of $X_*(A)$ by the coroot lattice for~$M$, and a homomorphism 
\[
w_M:M(L)
\to X_M.
\]  
For
$\mu,\nu \in X_M$ we write $\mu\overset{P}\le\nu$ if $\nu -\mu$ is a non-negative
integral linear combination of (images in~$X_M$ of) coroots $\alpha^\vee$, where
$\alpha$ ranges over the roots of
$A$ in~$N$.

We write
$A_P$ for the identity component of the center of~$M$; thus $A_P$ is a split torus
over~$F$. Let
$\mathfrak a_P$ denote the real vector space $X_*(A_P) \otimes_{\mathbf Z} \mathbf
R$. As usual
$P$ determines an open chamber $\mathfrak a_P^+$ in~$\mathfrak a_P$, defined by
\[
\mathfrak a_P^+:=\{ x \in \mathfrak a_P : \langle \alpha, x \rangle > 0 \text{ for
every root
$\alpha$ of $A_P$ in $N$} \}.
\] The composition $X_*(A_P) \hookrightarrow X_*(A) \twoheadrightarrow X_M$, when
tensored with
$\mathbf R$, yields a canonical isomorphism $\mathfrak a_P \simeq X_M
\otimes_\mathbf Z \mathbf R$. We write $X_M^+$ for the subset of~$X_M$ consisting
of all elements whose image in
$X_M\otimes_{\mathbf Z} \mathbf R \simeq \mathfrak a_P$ lies in~$\mathfrak a_P^+$. 

For any coweight $\mu \in X_*(A)$ (usually taken to be dominant) and any $b \in
G(L)$ we consider the affine Deligne-Lusztig set 
\begin{equation}\label{aDL}
 X^G_{\mu}(b):=\{ x \in G(L)/G(\mathfrak o_L) : x^{-1}b \sigma(x) \in G(\mathfrak
o_L)\mu(\varpi)G(\mathfrak o_L) \}.
\end{equation}
For $b \in M(L)$ we may replace  
 $G$ by $M$ in this definition, obtaining the affine
Deligne-Lusztig set $X^M_{\mu}(b)$; moreover, the inclusion $M(L)/M(\mathfrak o_L)
\hookrightarrow G(L)/G(\mathfrak o_L)$ induces an inclusion
\[
X^M_{\mu}(b)  \hookrightarrow X^G_{\mu}(b). 
\]

\begin{thm} \label{main.theorem}
 Let $\mu \in X_*(A)$ be a dominant coweight, and let $b$ be a basic
element in $M(L)$  such that
$w_M(b)$ lies in the subset
$X_M^+$ of~$X_M$.
\begin{enumerate}
\item \cite{rapoport-richartz96} If $X^G_{\mu}(b)$ is non-empty, then $w_M(b)
\overset{P} \le \mu$. Here we are regarding $\mu$ as an element of~$X_M$.
\item Suppose that $w_M(b)=\mu$ (with $\mu$ again regarded as an element
of~$X_M$). Then the natural injection $X^M_{\mu}(b) \hookrightarrow X^G_{\mu}(b)$
is a bijection.
\end{enumerate}
\end{thm}

See \cite{kottwitz85} for the definition of basic.
The first part of the theorem is a reformulation  of the
group-theoretic generalization of Mazur's inequality (see Theorem 1.4.1 of
\cite{katz79}) proved by Rapoport-Richartz \cite{rapoport-richartz96}. 
Proposition 
\ref{reform.prop} shows that our formulation 
 agrees with that of \cite{rapoport-richartz96}. The second part of
the theorem is the group-theoretic generalization of Katz's theorem (see Theorem
1.6.1 of
\cite{katz79}) which was formulated in
\cite[Remark 4.12]{kottwitz-rapoport02} and  proved  there for $GL_n$ and
$GSp_{2n}$ as a consequence of Katz's theorem. 

The theorem above is proved in
\S\ref{sec3}. In \S\ref{sec4} it is generalized to the case of unramified groups
(see Theorem \ref{thm4.2} for a precise statement). Theorem \ref{qsconv} 
determines (for any unramified $G$)  when the
affine Deligne-Lusztig set
$X^G_\mu(b)$ is non-empty, generalizing Proposition 4.6 of
\cite{kottwitz-rapoport02}, which treated the case of split groups.

We now draw the reader's attention to some related work.  In
\cite{fontaine-rapoport02} (respectively, \cite{kottwitz-rapoport02}) the converse
to Mazur's inequality is proved for $GL_n$ (respectively, $GL_n$ and $GSp_{2n}$).
Recently C.~Leigh
 \cite{leigh02} has  proved the converse to Mazur's inequality for all split
classical groups. The reader who would like to know how these questions relate to
the reduction modulo $p$ of Shimura varieties should consult the survey article
\cite{rapoport02} by Rapoport.

It is a pleasure to thank M.~Rapoport, both for the interest he has taken in this
work and  the helpful comments he made about a preliminary version of the paper.

\section{Retractions}\label{sec.retractions}
\subsection{Notation}\label{r.1} Let $\mathfrak o$ be a complete discrete
valuation ring with fraction field~$F$, uniformizing element $\varpi$, and residue
field $k=\mathfrak o/\varpi \mathfrak o$.

Let $G$ be a split connected reductive group over $\mathfrak o$ and let $A$ be a
split maximal torus of~$G$ over~$\mathfrak o$. We denote by $\mathcal B=\mathcal
B(A)$ the set of Borel subgroups of~$G$ containing $A$ (all of which are defined
over $\mathfrak o$). For $B \in \mathcal B$ denote by $\bar B$ the Borel subgroup
in $\mathcal B$ that is opposite to~$B$. For
$B$,$B_1$,$B_2,\dots$ in $\mathcal B$ we denote the unipotent radical by
$U$,$U_1$,$U_2,\dots$, so that (for instance) $B=AU$. We write $K$ for
$G(\mathfrak o)$.

\subsection{Definition of retractions} \label{r.2} For $g\in G(F)$ and $B=AU$
in~$\mathcal B$ we define $r_B(g) \in X_*(A)$ to be the unique element
$\mu \in X_*(A)$ such that $g \in K \cdot \mu(\varpi) \cdot U(F)$. The family
$(r_B(g))_{B \in
\mathcal B}$ of retractions is used by Arthur to form weighted orbital integrals.

\subsection{Positivity properties of families of retractions}\label{r.3} The
family $(r_B(g))_{B \in \mathcal B}$ has the following basic positivity property
\cite[Lemma 3.6]{arthur76},\cite[Lemma 85]{harish-chandra66b}. Let $B_1=AU_1$ and
$B_2=AU_2$ be adjacent Borel subgroups in~$\mathcal B$, and let $\alpha$ be the
unique root of~$A$ that is positive for~$B_1$ and negative for~$B_2$. Then 
\begin{equation}\label{r.3.1} r_{B_2}(g)-r_{B_1}(g)=j\cdot \alpha^\vee,
\end{equation} where $j$ is a non-negative integer that we will now define. (We
will not recall the proof except to say that one reduces to the case of $SL(2)$,
for which a simple computation with $2\times 2$ matrices does the job.) 

The group $U_1$ is the semidirect product of the normal subgroup $U_1 \cap U_2$
and the root subgroup $U_\alpha$ determined by $\alpha$. In particular $U_\alpha$
is a quotient of $U_1$, and we refer to the image of $u_1 \in U_1$ in $U_\alpha$
as the $\alpha$-component of $u_1$. Choosing an isomorphism between $U_\alpha$ and
$\mathbf G_a$ over $\mathfrak o$, we may view the
$\alpha$-component of $u_1 \in U_1(F)$ as an element of~$F$, well-defined up to
multiplication by a unit. 

Now we can define $j$. Decompose $g$ as $g=k\cdot u_1\cdot \mu(\varpi)$ with $k
\in K$, $u_1 \in U_1(F)$ and $\mu \in X_*(A)$ (so that $\mu=r_{B_1}(g)$), and
write $x \in F$ for the
$\alpha$-component of~$u_1$. Then $j$ is defined to be $0$ if $x \in \mathfrak o$ 
and is defined to be
$-\val(x)$ if $x \notin \mathfrak o$. 

The basic positivity property above has some obvious consequences. One is that for
any $B_1,B_2
\in \mathcal B$ the coweight $r_{B_2}(g)-r_{B_1}(g)$ is a non-negative integral
linear combination of coroots $\alpha^\vee$ that are positive for~$B_1$ and
negative for~$B_2$. Thus for any $B,B'
\in \mathcal B$ we have
\begin{equation}\label{r.3.2} r_{B}(g) \overset{B}\le r_{B'}(g) \overset{B}\le
r_{\bar B}(g),
\end{equation} where $\mu \overset{B}\le \nu$ (for coweights $\mu$,$\nu$) means
that $\nu-\mu$ is a non-negative integral linear combination of coroots that are
positive for~$B$. 

\subsection{Recognizing the subset $K\cdot M(F)$ of $G(F)$ using
retractions}\label{r.4} Let $M$ be a Levi subgroup of~$G$ containing $A$ and note
that $M$ is automatically defined over~$\mathfrak o$. We write $X_M$ for the
quotient of $X_*(A)$ by the coroot lattice for~$M$. For example, when $M=A$,
we have $X_M=X_*(A)$.

\begin{lem}\label{lemma.r.4.1} Let $g \in G(F)$. Then $g \in K \cdot M(F)$ if and
only if
$r_{B_1}(g)$ and $r_{B_2}(g)$ are equal in $X_M$ for all $B_1,B_2\in \mathcal B$.
\end{lem}
\begin{proof} Assume first that $g=km$ with $k \in K$ and $m \in M(F)$. We must
show that all the retractions of $g$ are equal in~$X_M$. For any $B \in \mathcal
B$ the intersection $B
\cap M$ is a Borel subgroup of~$M$, and it is clear that $r_B(g)=r_{B\cap M}(m)$.
Thus, for
$B_1,B_2 \in \mathcal B$ the coweight $r_{B_2}(g)-r_{B_1}(g)$ is a non-negative
integral linear combination of coroots $\alpha^\vee$ for~$M$ that are positive
for~$B_1$ and negative for~$B_2$, and in particular $r_{B_2}(g)=r_{B_1}(g)$ in
$X_M$.

Now assume that all the retractions of $g$ are equal in~$X_M$. Choose $B=AU$
in~$\mathcal B$ and choose a minimal gallery 
\begin{equation}
 B=B_0,B_1,\dots,B_{l-1},B_l=\bar B
\end{equation}
 of Borel subgroups in~$\mathcal B$ joining $B$ to~$\bar B$. Thus $l$
is equal to the number of positive roots for ~$B$, and the subgroups $B_i$,
$B_{i+1}$ are adjacent for $0 \le i \le l-1$. Write
$B_i=AU_i$ and put $V_i:=U \cap U_i$. Note that $V_0=U$ and $V_l=\{1\}$.
  We will prove by induction on $i$ (for $0 \le i \le l$) that $g \in K\cdot
V_i(F) \cdot M(F)$. The case $i=0$ is obvious and the case $i=l$ is the statement
of the lemma.

 For the induction step we suppose that for some $i$ less than $l$ we have $g=kum$
for $k\in K$, $u\in V_i(F)$, $m \in M(F)$. The group $V_i$ is the semidirect
product of the normal subgroup $V_{i+1}$ and the root subgroup
$U_\alpha$, where $\alpha$ is the unique root of $A$ that is positive for $B_i$
and negative for
$B_{i+1}$.  If $\alpha$ is not a root of~$M$, then
$\alpha^\vee$ is a non-torsion element in~$X_M$, and our hypothesis that the
retractions for $B_i$ and
$B_{i+1}$ are equal in $X_M$ ensures that the $\alpha$-component $u_\alpha$ of~$u$
lies in
$U_\alpha(\mathfrak o)$, so that we can write $g=(ku_\alpha)\cdot
u_\alpha^{-1}u\cdot m \in K\cdot V_{i+1}(F) \cdot M(F)$. If
$\alpha$ is a root of~$M$, then we can write $g =k\cdot uu_\alpha^{-1} \cdot
u_\alpha m\in K\cdot V_{i+1}(F)\cdot M(F)$.
\end{proof}

\subsection{Review of two relations between the Iwasawa and Cartan
decompositions}\label{r.6} We now recall two results of Bruhat-Tits. Fix $B \in
\mathcal B$ and a
$B$-dominant coweight $\mu \in X_*(A)$. Suppose that $g \in K\cdot \mu(\varpi)
\cdot K$. Then (see \cite[4.4.4]{bruhat-tits72}) for all $B' \in \mathcal B$
\begin{equation}\label{r.6.1} r_{B'}(g)\overset{B}\le \mu
\end{equation} and if $r_B(g)=\mu$, then $g \in K\cdot A(F)$.

Note that \eqref{r.3.2} and \eqref{r.6.1} together yield 
\begin{equation}\label{r.6.2} r_{B}(g) \overset{B}\le r_{B'}(g) \overset{B}\le \mu
\end{equation} for all $B' \in \mathcal B$, so that the hypothesis $r_B(g)=\mu$
implies that all the retractions of~$g$ are equal; therefore the second result of
Bruhat-Tits follows from the first together with Lemma \ref{lemma.r.4.1}. This
proof of their second result (different from the one given in
\cite{bruhat-tits72}) has the advantage that it generalizes immediately to
parabolic subgroups, as we now check.

\subsection{Variant (for parabolic subgroups) of the two results of
Bruhat-Tits}\label{r.7} Fix $B=AU \in \mathcal B$ as well as a parabolic
subgroup~$P$ of~$G$ containing~$B$. Write $P=MN$, where $M$ is the unique Levi
subgroup of~$P$ containing $A$ and $N$ is the unipotent radical of~$P$. As before
we write $X_M$ for the quotient of $X_*(A)$ by the coroot lattice for~$M$. For
$\mu,\nu \in X_M$ we write $\mu\overset{P}\le\nu$ if $\nu -\mu$ is a non-negative
integral linear combination of (images in~$X_M$ of) coroots $\alpha^\vee$, where
$\alpha$ ranges over the roots of
$A$ in~$N$ (or, equivalently, $U$).

As in the previous section, fix a
$B$-dominant coweight $\mu \in X_*(A)$ and an element $g \in K\cdot \mu(\varpi)
\cdot K$. It follows immediately from \eqref{r.6.1} that
\begin{equation}\label{r.7.1} r_{B'}(g)\overset{P}\le \mu
\end{equation} for all $B' \in \mathcal B$, where the two sides of this inequality
are now viewed as elements in~$X_M$.

\begin{lem}\label{lemma.r.7.1} With $\mu$, $g$ as above assume further that
$r_B(g)$ is equal to $\mu$ in~$X_M$. Then $g \in K
\cdot M(F)$. Moreover, writing $g=km$ for $k \in K$ and $m \in M(F)$, then $m \in
K_M \cdot
\mu(\varpi) \cdot K_M$, where we have written $K_M$ for $M(\mathfrak o)$.
\end{lem}
\begin{proof} Let $B' \in \mathcal B$. Since (by \eqref{r.6.2})
\begin{equation}\label{r.7.2} r_{B}(g) \overset{P}\le r_{B'}(g) \overset{P}\le \mu,
\end{equation} our hypothesis that $r_B(g)=\mu$ in~$X_M$ implies that
$r_B(g)=r_{B'}(g)$ in~$X_M$ for all $B' \in
\mathcal B$. Therefore Lemma \ref{lemma.r.4.1} implies that $g \in K \cdot M(F)$. 

Now we write $g$ as $km$ and verify the second statement of the lemma. By the
Cartan decomposition for~$M$ there exists a unique coweight $\nu \in X_*(A)$ that
is dominant with respect to the Borel subgroup $B \cap M$ of~$M$ and is such that
$m$ lies in $K_M \cdot \nu(\varpi) \cdot K_M$. By the Cartan decomposition for~$G$
the coweights $\mu$ and $\nu$ lie in the same orbit of the Weyl group of~$A$
in~$G$. Since both $\mu$ and $\nu$ are dominant for~$M$, no root hyperplane
for~$M$ separates $\mu$ from~$\nu$. Therefore $\mu - \nu$ is a non-negative
integral linear combination of coroots
$\alpha^\vee$, for $\alpha$ ranging through the roots of~$A$ in~$N$; on the other
hand it is clear that $\mu$ and $\nu$ are equal in~$X_M$ (since
$\mu=r_B(g)=r_{B\cap M}(m)$ and $m \in K_M\nu(\varpi)K_M)$; therefore $\mu=\nu$,
showing that $m \in K_M \cdot
\mu(\varpi) \cdot K_M$, as desired.
\end{proof}

\section{Generalizations of Mazur's inequality and Katz's theorem}\label{sec3}
\subsection{Notation}\label{m.1} In the rest of the paper $F$
denotes a finite extension of~$\mathbf Q_p$ and $\mathfrak o$ denotes the
valuation ring of~$F$. We write
$L$ for the completion of the maximal unramified extension of~$F$ in some
algebraic closure~$\bar F$ of~$F$. We write $\sigma$ for the Frobenius
automorphism of~$L$ over~$F$, and we write 
$\mathfrak o_L$ for the valuation ring of~$L$. 

\subsection{A lemma about $\sigma$-$L$-spaces} \label{m.2} Recall that a
$\sigma$-$L$-\emph{space} is a pair $(V,\Phi)$ consisting of a finite dimensional
vector space $V$ over~$L$ and a $\sigma$-semilinear bijection $\Phi:V \to V$. In
case $F=\mathbf Q_p$ a $\sigma$-$L$-space is an isocrystal, and the theory of
$\sigma$-$L$-spaces is completely parallel to that of isocrystals. In particular
there are finitely many rational numbers, called slopes, attached to $(V,\Phi)$
(see \S3 in \cite{kottwitz85}). 
\begin{lem}\label{lemma.m.2.1}  Let $(V,\Phi)$ be a $\sigma$-$L$-space and assume
that all its slopes are strictly positive. 
\begin{enumerate}
\item For any $v \in V$ the sequence $\Phi^n v$ approaches $0$ as $n \to +\infty$.
\item Suppose that $\Lambda$ is an $\mathfrak o_L$-lattice in~$V$ such that $\Phi
\Lambda \subset
\Lambda$, and suppose that $v$ is an element of~$V$ such that $v-\Phi v \in
\Lambda$. Then $v \in
\Lambda$. 
\end{enumerate}
\end{lem}
\begin{proof} We begin by proving the first part of the lemma. Choose a positive
integer $j$ such that $jr \in
\mathbf Z$ for every slope $r$ of $(V,\Phi)$. Then, in a suitable basis for~$V$
the map $\Phi^j$ can be represented by a diagonal matrix whose diagonal entries
are strictly positive powers of the uniformizing parameter $\varpi$ for ~$F$, and
it is clear that $\Phi^{jm}v' \to 0$ as $m \to +\infty$ for every $v' \in V$.
Taking for $v'$ the $j$ vectors $v,\Phi v,\dots,\Phi^{j-1}v$, we see that $\Phi^n
v \to 0$, as desired.  Now we prove the second part of the lemma. It follows from
the first part of the lemma that we may define an additive homomorphism $\Psi:V
\to V$ by $\Psi=1+\Phi + \Phi^2 + \Phi^3 + \dots$ and hence  that the additive
homomorphism $1-\Phi$ is bijective with inverse $\Psi$. Also, it is clear from the
definition of $\Psi$ that $\Psi \Lambda \subset \Lambda$. 

We are given $v \in V$ such that $(1-\Phi)v \in \Lambda$. Applying $\Psi$, we
conclude that $v \in
\Psi \Lambda \subset \Lambda$, as desired. 
\end{proof}

\subsection{Proof of Theorem
\ref{main.theorem}}\label{m.3}    In the proof of Theorem
\ref{main.theorem} we will need the following non-abelian analog of  Lemma
\ref{lemma.m.2.1}.

\begin{lem}\label{lemma.m.3.1} Let $\mu \in X_*(A)$
be a dominant coweight, and let $b$ be a basic element in $M(L)$ 
 such that $w_M(b)$ lies in the
subset
$X_M^+$ of~$X_M$. Assume further that 
$b \in M(\mathfrak o_L)
\mu(\varpi)M(\mathfrak o_L) $. Write $\Phi$ for the automorphism $n \mapsto
b\sigma(n)b^{-1}$ of~$N$ over~$L$. Let $n \in N(L)$ and assume that $n^{-1}\Phi(n)
\in N(\mathfrak o_L)$. Then $n \in N(\mathfrak o_L)$. 
\end{lem}
\begin{proof} Conjugation by $M(\mathfrak o_L)$ preserves $N(\mathfrak o_L)$, and
since $\mu$ is dominant, we have $\mu(\varpi)N(\mathfrak o_L) \mu(\varpi)^{-1}
\subset N(\mathfrak o_L)$; it follows that
$bN(\mathfrak o_L)b^{-1}
\subset N(\mathfrak o_L)$ and hence that $\Phi N(\mathfrak o_L) \subset
N(\mathfrak o_L)$. Our hypothesis that $w_M(b) \in X_M^+$ ensures that all the
slopes of $\Phi$ on $\Lie N(L)$ are strictly positive. 
Indeed, since $b$ is basic in~$M(L)$, these slopes are given by $\langle \alpha,
\overline{w_M(b)}
\rangle$, where $\alpha$ ranges through the roots of~$A_P$ on~$\Lie(N)$, and
$\overline{w_M(b)}$ denotes the image of $w_M(b)$ in~$\mathfrak a_P$.  Thus, if $N$
is abelian, we have only to appeal to Lemma
\ref{lemma.m.2.1}.

In the general case we need to choose an $M$-stable filtration 
\begin{equation}
 N=N_0 \supset N_1 \supset N_2 \supset \dots \supset N_r=\{1\}
\end{equation}
 by normal subgroups with $N_i/N_{i+1}$ abelian for all $i$. Each
$N_i$ is
$A$-stable, hence is a product of root subgroups and is necessarily defined over
$\mathfrak o$. We will prove by induction on $i$ ($0 \le i \le r$) that $n \in
N_i(L)\cdot N(\mathfrak o_L)$. For $i=0$ this statement is trivial, and for $i=r$
it is the statement of the lemma. It remains to do the induction step. So suppose
that for $0 \le i  < r$ we can write $n$ as $n=n_i n_{\mathfrak o}$ for
$n_i \in N_i(L)$ and $n_\mathfrak o \in N(\mathfrak o_L)$. Then $n_i^{-1}\Phi(n_i)
\in N_i(\mathfrak o_L)$. So (by Lemma \ref{lemma.m.2.1}) the image of $n_i$ in
$(N_i/N_{i+1})(L)$ lies in $(N_i/N_{i+1})(\mathfrak o_L)$. Since $N_i(\mathfrak
o_L)$ maps onto $(N_i/N_{i+1})(\mathfrak o_L)$, we see that $n_i$ can be written
as $n_{i+1} n'_{\mathfrak o}$ with $n_{i+1} \in N_{i+1}(L)$ and $n'_{\mathfrak o}
\in N_i(\mathfrak o_L)$. Thus $n=n_{i+1} \cdot (n'_\mathfrak o  n_\mathfrak o) \in
N_{i+1}(L)N(\mathfrak o_L)$, as desired.
\end{proof}

Now we are ready to prove the main theorem for split groups.
\begin{proof}[Proof of Theorem \ref{main.theorem}] Let $g \in G(L)$ and suppose
that 
\begin{equation}\label{m.3.1} g^{-1}b\sigma(g) \in K_L\mu(\varpi)K_L,
\end{equation}  where we have written
$K_L$ for~$G(\mathfrak o_L)$. Use the Iwasawa decomposition to write $g$ as $mnk$
for $m \in M(L)$, $n \in N(L)$ and $k \in K_L$. It follows from \eqref{m.3.1} that
\begin{equation}\label{m.3.2} n_1m_1 \in K_L\mu(\varpi)K_L,
\end{equation} where $m_1:=m^{-1}b\sigma(m) \in M(L)$ and
$n_1:=n^{-1}m_1\sigma(n)m_1^{-1} \in N(L)$. We claim that 
\begin{equation}\label{353}
w_M(b)=r_B(n_1m_1),
\end{equation}
with the right side being regarded as an element of~$X_M$. Indeed,
\begin{equation} w_M(b)=w_M(m_1)=r_{B\cap M}(m_1),
\end{equation} which in turn is equal to the image in~$X_M$ of $r_B(n_1m_1)$. 

We
conclude from \eqref{r.7.1}, \eqref{m.3.2}, \eqref{353} that
$w_M(b) \overset{P} \le \mu$, which proves the first part of the theorem. 

Now we prove the second part of the theorem. Under the
hypothesis that $w_M(b)=\mu$ (and with 
$g,m,n,m_1,n_1$ as above), we must prove that $g \in M(L) \cdot K_L$. It
follows from \eqref{m.3.2}, \eqref{353} and Lemma
\ref{lemma.r.7.1} that
$n_1m_1 \in K_L \cdot M(L)$. Therefore $n_1 \in K_L\cdot M(L)$, say $n_1=k_2m_2$
with $k_2 \in K_L$ and $m_2 \in M(L)$. Then $n_1m_2^{-1} \in P(\mathfrak o_L)$,
and therefore $n_1 \in N(\mathfrak o_L)$ and $m_2 \in M(\mathfrak o_L)$. Since
$n_1 \in N(\mathfrak o_L)$, the second statement of Lemma \ref{lemma.r.7.1}
applies to $n_1m_1$, and hence $m_1 \in M(\mathfrak
o_L)\mu(\varpi)M(\mathfrak o_L)$. 

Now applying Lemma \ref{lemma.m.3.1} (not to the element $b$, but to its
$\sigma$-conjugate $m_1$, which satisfies the same hypotheses as $b$), we see that
$n \in N(\mathfrak o_L)$. Therefore $g= m \cdot nk \in M(L)K_L$, and we are done,
since we have already seen that 
\begin{equation} m^{-1}b\sigma(m)=m_1
\in M(\mathfrak o_L)\mu(\varpi)M(\mathfrak o_L).
\end{equation}
\end{proof}

\section{Unramified groups}\label{sec4} It is easy to generalize  Theorem
\ref{main.theorem} from the case of split groups to that of unramified groups,
in other words, quasi-split groups $G$ over~$F$ that split over $L$. There is no
need to generalize the results in \S\ref{sec.retractions}, since we will apply
them to the group $G$ over~$L$, where it becomes split. 

Continuing with unramified $G$, we will then determine precisely which affine
Deligne-Lusztig sets \eqref{aDL} are non-empty, generalizing Proposition 4.6 of
\cite{kottwitz-rapoport02}.  

\subsection{Notation}\label{not} 
We will now change notation slightly, to emphasize
that our maximal torus is no longer assumed to be split. We consider a Borel
subgroup
$B=TU$ of~$G$, where
$T$ is a maximal torus in~$B$ and $U$ is the unipotent radical of~$B$; all these
subgroups are assumed to be defined over~$\mathfrak o$. In addition we fix a
parabolic subgroup $P=MN$ containing $B$, with $M$ containing $T$; again 
all these subgroups are assumed to be defined over~$\mathfrak o$. 

We denote by $A_P$ the maximal split torus in the center of~$M$, and we write
$\mathfrak a_P$ for $X_*(A_P) \otimes_{\mathbf Z} \mathbf R$.  In the special case
$P=B$, we often write $A$ and $\mathfrak a$ rather than $A_B$ and $\mathfrak a_B$;
of course $A$ is simply the maximal split torus in~$T$. As usual we identify
$\mathfrak a_P$ with a subspace of $\mathfrak a$. 

As before we denote by $X_M$ the quotient of $X_*(T)$ by the coroot lattice
for~$M$. The Frobenius automorphism $\sigma$ acts on $X_M$, and we denote by $Y_M$
the coinvariants of~$\sigma$ on~$X_M$. Thus $Y_M$ is the quotient of~$X_M$ by the
image of the homomorphism $1-\sigma:X_M \to X_M$. We introduce a partial order
on $Y_M$ as follows: for $y_1,y_2 \in Y_M$ we say that $y_2 \overset{P}
\preccurlyeq y_1$ if $y_1-y_2$ is a non-negative integral linear combination of
images in~$Y_M$ of coroots $\alpha^\vee$ corresponding to simple roots $\alpha$
of~$T$ in~$N$. 

We identify $Y_M
\otimes_{\mathbf Z} \mathbf R$ with $\mathfrak a_P$, and we write $Y^+_M$ for the
subset of $Y_M$ consisting of those elements whose image in $\mathfrak a_P$ lies
in the cone 
\begin{equation}
\mathfrak a^+_P:=\{x \in \mathfrak a_P : \langle \alpha,x \rangle > 0
\text{ $\forall$  root $\alpha$ of $A_P$ in $N$ } \}.
\end{equation}

As in \cite{kottwitz85}) the homomorphism $w_M:M(L) \to X_M$ induces a map
\begin{equation}
\kappa_M:B(M) \to Y_M.
\end{equation}
Now we can generalize Theorem \ref{main.theorem} to this more general context,
with the affine Deligne-Lusztig sets $X^G_\mu(b)$ still defined by \eqref{aDL}.

\begin{thm}\label{thm4.2}Let $\mu \in X_*(T)$ be a dominant coweight, and let $b$
be a basic element in
$M(L)$ such that $\kappa_M(b) \in Y_M^+$. 
\begin{enumerate}
\item \cite{rapoport-richartz96} \label{qs1}
 If $X^G_{\mu}(b)$ is non-empty, then $\kappa_M(b)
\overset{P} \preccurlyeq \mu$. Here we are regarding $\mu$ as an element of~$Y_M$.
\item \label{qs2}
Suppose that $\kappa_M(b)=\mu$ (with $\mu$ again regarded as an element
of~$Y_M$). Then the natural injection $X^M_{\mu}(b) \hookrightarrow X^G_{\mu}(b)$
is a bijection.
\end{enumerate}
\end{thm}
\begin{proof}
The proof of part (\ref{qs1}) is the same as the proof in the split case. The proof
of part (\ref{qs2}) is the same as in the split case, except for one new
complication, which we will now explain. 
 Suppose that $g \in G(L)$ represents an element of $X^G_\mu(b)$, so 
that 
\begin{equation}\label{q.3.1} g^{-1}b\sigma(g) \in K_L\mu(\varpi)K_L,
\end{equation}  where we have written
$K_L$ for~$G(\mathfrak o_L)$. As before we use the Iwasawa decomposition to write
$g$ as
$mnk$ for $m \in M(L)$, $n \in N(L)$ and $k \in K_L$. It follows from
\eqref{q.3.1} that
\begin{equation}\label{q.3.2} n_1m_1 \in K_L\mu(\varpi)K_L,
\end{equation} where $m_1:=m^{-1}b\sigma(m) \in M(L)$ and
$n_1:=n^{-1}m_1\sigma(n)m_1^{-1} \in N(L)$. We need to prove that
$g \in M(L)
\cdot K_L$ (under the hypothesis that $\kappa_M(b)=\mu$). 
Denote by $\nu$ the image of $r_B(n_1m_1)$ in~$X_M$. 
As in the proof of
Theorem \ref{main.theorem}, the elements $\nu$, $w_M(b)$ of $X_M$ have the same
image in~$Y_M$, and by hypothesis the image of $w_M(b)$ in $Y_M$ is $\mu$. We would
like to apply Lemma \ref{lemma.r.7.1} to the element $n_1m_1$, but for this we
would need to know that $\nu$ and $\mu$ are equal in~$X_M$, while all we know at
the moment is that they are equal in the quotient $Y_M$ of~$X_M$. However, by 
\eqref{q.3.2} and \eqref{r.7.1} we also know that  $ \nu \overset{P} \le \mu$.
Therefore Lemma
\ref{complication} shows that $\nu=\mu$ in $X_M$, as desired. Thus we may apply
Lemma \ref{lemma.r.7.1} to the element $n_1m_1$ in order to see that $n_1m_1 \in
K_L \cdot M(L)$. The rest of the proof is exactly the same as in the split case.
Of course we need to appeal to Lemma \ref{lemma.m.3.1}, but its statement
generalizes without change to the general unramified case, and its proof stays the
same too, though one should note that the subgroups $N_i$ used in the proof need
to be  chosen so as to be defined over 
$\mathfrak o$. 
\end{proof}
Next we prove the lemma we just used.

\begin{lem}\label{complication}
Let $x \in X_M$. Suppose that $x \overset{P} \ge 0$ and suppose further that the
image of~$x$ in~$Y_M$ is $0$. Then $x=0$. 
\end{lem}
\begin{proof}
Let $X^G_M$ denote the kernel of the canonical surjection $X_M \to X_G$, so that 
we get a short exact sequence
\begin{equation}
0 \to X^G_M \to X_M \to X_G \to 0.
\end{equation}
Taking coinvariants for $\sigma$, we get an exact sequence
\begin{equation}
 Y^G_M \to Y_M \to Y_G \to 0,
\end{equation}
where $Y^G_M$ denotes the coinvariants of $\sigma$ on~$X^G_M$. 

Clearly $X^G_M$ is a free abelian group on the set $S$ of coroots $\alpha^\vee$
for simple roots $\alpha$ of~$T$ that occur in~$N$, and $\sigma$ permutes these
basis elements. Therefore $Y^G_M$ is a free abelian group on the set $\bar S$ of
orbits of $\sigma$ on~$S$. In particular $Y^G_M$ is torsion-free, which implies
that the map $Y^G_M \to Y_M$ is injective. 

Now consider $x \in X_M$ such that $x \overset{P} \ge 0$. In particular $x$ lies
in the subgroup $X^G_M$.  Consider the image $y$ of $x$ in~$Y_M$. It is clear that
$y$ lies in the subgroup $Y^G_M$, and that its coefficients in the basis $\bar S$
are given by summing over the orbits of~$\sigma$ on~$S$ the coefficients of~$x$ in
the basis~$S$; since these latter coefficients are non-negative by our hypothesis
that $x \overset{P} \ge 0$, we see that if $y=0$, then $x=0$ as well.
\end{proof}

\subsection{Non-emptiness of affine Deligne-Lusztig subsets of $G(L)/G(\mathfrak
o_L)$} We will now determine when the affine Deligne-Lusztig set $X^G_\mu(b)$ (see
\eqref{aDL}) is non-empty. The split case is treated in Proposition 4.6 of
\cite{kottwitz-rapoport02}, and we follow closely the method used there.

As in \cite{kottwitz-rapoport02}, for $\mu \in X_*(T)$ we put
\begin{equation}
\mathcal P_\mu:=\{ \nu \in X_*(T): \nu=\mu \text{ in $X_G$, and $\nu \in
\Conv(W\mu)$} \}
\end{equation}
where $\Conv(W\mu)$ denotes the convex hull of the $W$-orbit $W\mu$ of~$\mu$ in
$X_*(T)\otimes \mathbf R$, $W$ being the absolute Weyl group of~$T$ in~$G$.
 With $P=MN$, $Y_M$ as in \ref{not} we let $\mathcal P_{\mu,M}$ denote the subset
of
$Y_M$ obtained as the image of $\mathcal P_\mu$ under the canonical surjection
$X_*(T) \to Y_M$. 

\begin{thm}\label{qsconv}
Let $b \in M(L)$ be basic, and let $\mu \in X_*(T)$.  Then $X^G_\mu(b)$ is
non-empty if and only if  $\kappa_M(b)$ lies in  $\mathcal
P_{\mu,M}$. 
\end{thm}

Since every $\sigma$-conjugacy class in $G(L)$ contains an element that is basic in
some standard Levi subgroup~$M$ (see \cite{kottwitz85}), this theorem gives a
complete answer to the question of when $X^G_\mu(b)$ is non-empty. 
The theorem follows immediately from the following three lemmas.

\begin{lem}
Let $b \in G(L)$ be basic. Then the $\sigma$-conjugacy class of~$b$ meets $K_L
\mu(\varpi) K_L$ if and only if $\kappa_G(b)$ is equal to the image of~$\mu$
in~$Y_M$. 
\end{lem}
\begin{proof}
($\Longrightarrow$) Obvious, since the homomorphism $w_G$ is trivial on $K_L$.

($\Longleftarrow$) We write $N(T)$ for the normalizer of~$T$ in~$G$. Then we have
an exact sequence
\begin{equation}
1 \to T(\mathfrak o_L) \to N(T)(L) \to \tilde W \to 1,
\end{equation}
where $\tilde W:=W \ltimes X_*(T)$ is the extended affine Weyl group for~$G$
over~$L$. Using this exact sequence, one sees without difficulty that the map
$N(T)(L) \to \tilde W$ induces a bijection from $B(N(T))$ to the set of
$\sigma$-conjugacy classes in the affine Weyl group $\tilde W$. 

For $w \in W$ choose $\dot w \in N(T)(\mathfrak o_L)$ such that $\dot w \mapsto
w$. Associated (see \cite{kottwitz85}) to the element $b':=\mu (\varpi)\dot w$ 
of 
$N(T)(L)$ is a homomorphism $\nu: \mathbf D \to N(T)$, where $\mathbf D$ is the
diagonalizable group with character group $\mathbf Q$. It is easy 
to calculate $\nu$. Indeed, choose a positive integer $r$ such that
$(w\sigma)^r=\sigma^r$. We use $w$ to twist the action of $\sigma$ on~$T$,
obtaining a new unramified torus $T_w$ which becomes equal to $T$ over the
fixed field $F_r$ of $\sigma^r$ in~$L$, but for which the action of $\sigma$ is
now given by $w\sigma$ rather than $\sigma$. Thus $X_*(T_w)$ coincides with
$X_*(T)$ as an abelian group, but $\sigma$ acts by $w\sigma$ rather than $\sigma$.
From this point of view, the homomorphism $\nu:\mathbf D \to N(T)$ may be regarded
as an element of $\Hom(\mathbf D,T_w)$, or in other words as an 
 element in $X_*(T)\otimes \mathbf Q$  fixed by $w\sigma$, the explicit formula
for $\nu$ being given by 
\begin{equation}
r^{-1}\sum_{i=1}^r (w\sigma)^i(\mu) \in X_*(T)\otimes \mathbf Q.
\end{equation} 

The remarks above are valid for any $w \in W$. The well-known fact that
semisimple groups over finite fields have anisotropic maximal tori translates 
into the fact that there exists
$w \in W$ for which the torus $T_w$ is anisotropic modulo the center of~$G$; for
the rest of this proof we work with such an element $w$. In this case $\nu$
is forced to be central in~$G$, and therefore the element $b'=\mu (\varpi)\dot w$
is basic in~$G(L)$. It is obvious from the form of $b'$ that $\kappa_G(b')=\mu$.
By hypothesis $b$ is basic and $\kappa_G(b)=\mu$. Therefore \cite{kottwitz85}
$b$ is $\sigma$-conjugate in~$G(L)$ to $b'$. From the form of $b'$ it is obvious
that $b'$ lies in the $K_L$-double coset of $\mu(\varpi)$. This concludes the
proof. 
\end{proof}

As in \cite{kottwitz-rapoport02} for $\mu \in X_*(T)$ we denote by $M(\mu)$ the
image of $K_L \mu(\varpi)K_L \cap P(L)$ under the canonical surjection $P(L) \to
M(L)$; obviously $M(\mu)$ is a union of $M(\mathfrak o_L)$-double cosets. 

\begin{lem}
Let $b \in M(L)$ and let $\mu \in X_*(T)$. Then the following three conditions are
equivalent: 
\begin{enumerate}
\item The $\sigma$-conjugacy class of $b$ in $G(L)$ meets $K_L\mu(\varpi)K_L$.
\item The $\sigma$-conjugacy class of $b$ in $P(L)$ meets $K_L \mu(\varpi)K_L \cap
P(L)$.
\item The $\sigma$-conjugacy class of $b$ in $M(L)$ meets $M(\mu)$.
\end{enumerate}
\end{lem}
\begin{proof}
The equivalence of the first two conditions is clear from the Iwasawa
decomposition $G(L)=K_LP(L)$. The equivalence of the second two conditions follows
from the fact \cite[3.6]{kottwitz97} that $B(P) \to B(M)$ is a bijection.
\end{proof}

\begin{lem}
The set  $w_M(M(\mu))$ is equal to the image of $\mathcal P_\mu$ under the
canonical surjection  $X_*(T) \to X_M$.
\end{lem}
\begin{proof}
This is Lemma 4.5 of \cite{kottwitz-rapoport02}, which can be applied since $G$
splits over $L$. Note that the assumption, made at the beginning of \S4 of
\cite{kottwitz-rapoport02}, that the derived group of
$G$ be simply connected was made merely for convenience; in particular it was not
used in the proof of Lemma 4.5.
\end{proof}

\subsection{Remarks concerning the converse to Mazur's inequality}
Let $b$, $\mu$ be as in Theorem \ref{thm4.2}. Mazur's inequality (in other words
the first part of that theorem) says that if $X^G_{\mu}(b)$ is non-empty, then
$\kappa_M(b)
\overset{P} \preccurlyeq \mu$, where $\mu$ is being regarded an element of~$Y_M$. 
 Thus the converse to Mazur's inequality is the statement, only
known to be true in certain cases, that if $\kappa_M(b)
\overset{P} \preccurlyeq \mu$, then $X^G_{\mu}(b)$ is non-empty. Since Theorem
\ref{qsconv} does tell us exactly when $X^G_\mu(b)$ is non-empty, proving the
converse to Mazur's inequality is the same as proving that 
$\kappa_M(b)
\overset{P} \preccurlyeq \mu$ 
is equivalent to 
$\kappa_M(b) \in \mathcal
P_{\mu,M}$.
Thus, in order to prove the converse to Mazur's inequality in general, it would be
enough to answer the following question about root systems affirmatively.

\begin{question} \label{quest}
Let $\mu \in X_*(T)$ be a dominant coweight and let $\nu \in Y_M^+$. Are the
following two conditions equivalent?
\begin{enumerate} 
\item $\nu \overset{P}\preccurlyeq \mu$ \label{cond1}
\item $\nu \in \mathcal P_{\mu,M}$
\end{enumerate}
\end{question}
It is immediate that the second condition implies the first. The challenge is to
prove that the first condition implies the second; this was done for $GL_n$ and
$GSp_{2n}$ in \cite{kottwitz-rapoport02} and for all split classical groups in
\cite{leigh02};  non-split groups have not been examined yet.

\subsection{Comparison with \cite[4.2]{rapoport-richartz96}}\label{reform}
Condition \eqref{cond1} in Question \ref{quest} looks superficially different from
the one used in \cite{rapoport-richartz96} (and also used to define the set
$B(G,\mu)$ in \cite{kottwitz97}), but in fact it is equivalent, as we now check 
(see Proposition \ref{reform.prop} below). 

This equivalent condition involves some additional notation. We write
$X_\mathbf R$ for the real vector space $X_*(T) \otimes \mathbf R$. We identify 
$\mathfrak a$ with the subspace of $\sigma$-fixed vectors in $X_\mathbf R$, and we
view $\mathfrak a$ as a direct summand of $X_\mathbf R$,  the projection map
$X_ \mathbf R \to \mathfrak a$, denoted $x \mapsto x^\flat$, being given by
averaging over orbits of~$\sigma$. We have already identified $\mathfrak a_P$ with
a subspace of $\mathfrak a$. In fact we view $\mathfrak a_P$ as a direct summand of
$\mathfrak a$, the projection map $pr_M:\mathfrak a \to \mathfrak a_P$  being
given by averaging over the relative Weyl group $W_{M(F)}$ of~$T$ in~$M$. As usual
we identify $W_{M(F)}$ with the fixed points of $\sigma$ in the absolute Weyl
group $W_M$ of~$T$ in~$M$.

 The partial order $\overset{B}\le$ on~$X_*(T)$
extends as usual to a partial order on~$X_\mathbf R$, which we will denote simply
by $\le$; thus for
$x,y \in
X_\mathbf R$ the inequality $x \le y$ means that $y-x$ is a non-negative real
linear combination of simple coroots. 

Recall that we have already identified $Y_M \otimes \mathbf R$ with $\mathfrak
a_P$. For $\nu \in Y_M$ we denote by $\bar \nu$ the image of $\nu$ in $\mathfrak
a_P$ (which lets us view $\bar \nu$ as an element in $X_\mathbf R$, as we will do
in the next lemma). 

\begin{lem}\label{po}
Let $\nu \in Y_M$. Then $\nu \overset{P} \succcurlyeq 0 $ if and only if $\bar \nu
\ge 0$ and the image of $\nu$ in $Y_G$ is $0$.
\begin{proof}
Exercise.
\end{proof}
\end{lem}

\begin{lem}\label{dom} Let $x$ be a dominant element of~$X_\mathbf R$, and 
let $y \in \mathfrak a_P$. Then $y \le x^\flat$ if and only if
   $y \le pr_M(x^\flat) $.
\begin{proof} ($\Longrightarrow$) Apply $pr_M$ to the inequality $y \le x^\flat$, 
using that  $pr_M$ preserves $\le$.

($\Longleftarrow$)
This follows from the fact that $pr_M(x^\flat) \le x^\flat$, a consequence of the
inequalities  $wx^\flat \le x^\flat$ for  $w \in W_{M(F)}$, which hold since
$x^\flat $ is dominant.
\end{proof}
\end{lem}

\begin{prop} \label{reform.prop}
Let $\mu \in X_*(T)$ be a dominant coweight and let $\nu \in Y_M$. Then the 
following two conditions are equivalent:
\begin{enumerate}
\item $\nu \overset{P}\preccurlyeq \mu$. \label{l1}
\item \label{l2}
$\mu,\nu$ have the same image in $Y_G$, and  $\bar \nu \le \mu^\flat$. 
\end{enumerate}
\end{prop}
\begin{proof}
This follows from Lemmas \ref{po} and \ref{dom}.
\end{proof}

\bibliographystyle{amsalpha}
\providecommand{\bysame}{\leavevmode\hbox to3em{\hrulefill}\thinspace}

\end{document}